\documentclass[12pt]{article}
 
\usepackage[margin=1in]{geometry} 
\usepackage{amsmath,amsthm,amssymb}

\usepackage[latin1]{inputenc}
\usepackage{color} 
\usepackage{indentfirst}

\newcommand{\dis}{\displaystyle}
\newcommand{\R}{\mathbb{R}}
\newcommand{\C}{\mathbb{C}}

\newcommand{\N}{\mathbb{N}}

\def\res{\mathop{\operatorname{res}}}

\newcommand{\M}{\mathcal{M}}

\newcommand{\bigO}{\mathcal{O}}

\newcommand{\Real}{\operatorname{Re}}
\newcommand{\Imag}{\operatorname{Im}}

\newtheorem{theorem}{Theorem}[section]
\newtheorem{proposition}[theorem]{Proposition}

\newtheorem{definition}[theorem]{Definition}
\newtheorem{remark}[theorem]{Remark}

\numberwithin{equation}{section}
\begin{document}
\title{On Müntz-type formulas related to the Riemann zeta function}
\author{Hélder Lima}
\date{}
\maketitle 

\uppercase{Abstract.}
The Mellin transform and several Dirichlet series related with the Riemann zeta function are used to deduce some identities similar to the classical Müntz formula \cite{TitchmarshRiemann}.
These formulas are derived in the critical strip and in the half-plane $\Real(s)<0$.
As particular cases, integral representations for products of the gamma and zeta functions are exhibited.

\textbf{Keywords:} \textit{Arithmetic functions, Dirichlet series, Mellin transform, Riemann zeta function, Euler gamma function, Müntz formula, Müntz-type formulas.}

\textbf{AMS Subject Classifications:} 11M06, 11M26, 33B15, 42A38, 42B10, 44A05.   
 
\newpage
\section{Introduction}
The Mellin transform \cite{TitchmarshFourier} of a function $f$ is defined by
\begin{align}
\label{Mellin transform definition}
f^*(s)=\int_{0}^{\infty}f(x)x^{s-1}dx
\end{align}
and its inverse transform is given by 
\begin{align}
\label{Mellin inverse transform}
f(x)=\frac{1}{2\pi i}
\int_{\sigma-i\infty}^{\sigma+i\infty}f^*(s)x^{-s}ds.
\end{align}

The following proposition establishes reciprocity between (\ref{Mellin transform definition}) and (\ref{Mellin inverse transform}) under favourable conditions (see section 1.29 of \cite{TitchmarshFourier}).

\begin{proposition}
\label{Mellin transform of a Mellin inverse transform}
Suppose that $a,b\in\R$, $a<b$, $\dis f^*(s)$ is an analytic function in the strip $a<\Real(s)<b$ such that, for each $a<\sigma<b$, $\dis f^*(s)\in L_1(\sigma\pm i\infty)\equiv L_1(\sigma- i\infty,\sigma+ i\infty)$ and $f(x)$ is defined by (\ref{Mellin inverse transform}).
Then $\dis f(x)x^{\sigma-1}\in L_1(0,\infty)$, for all $a<\sigma<b$, and the Mellin transform of $f$ is equal to $\dis f^*(s)$ in the strip $a<\Real(s)<b$.
\end{proposition}

As it is known, the Riemann zeta function ($\zeta(s)$) \cite{TitchmarshRiemann} is analytic in the entire complex plane except the point $s=1$, where it has a simple pole such that $\dis\res_{s=1}\zeta(s)=\lim_{s\to 1}(s-1)\zeta(s)=1$.
Moreover, we have the following representations of expressions involving the Riemann zeta function in form of Dirichlet series absolutely convergent in the half-plane $\Real(s)>1$:
\begin{align}
\label{zeta function definition}
\zeta(s)=\sum_{n=1}^{\infty}\frac{1}{n^s};
\end{align}

\begin{align}
\label{Dirichlet 1/zeta}
\frac{1}{\zeta(s)}=\sum_{n=1}^{\infty}\frac{\mu(n)}{n^s};
\end{align}

\begin{align}
\label{Dirichlet zeta^k(s)}
\zeta^k(s)=\sum_{n=1}^{\infty}\frac{d_k(n)}{n^s},\,k\in\N;
\end{align}

\begin{align}
\label{Dirichlet zeta(s)/zeta(2s)}
\frac{\zeta(s)}{\zeta(2s)}
=\sum_{n=1}^{\infty}\frac{|\mu(n)|}{n^s};
\end{align}

\begin{align}
\label{Dirichlet zeta^2(s)/zeta(2s)}
\frac{\zeta^2(s)}{\zeta(2s)}
=\sum_{n=1}^{\infty}\frac{2^{\omega(n)}}{n^s};
\end{align}

\begin{align}
\label{Dirichlet zeta^3(s)/zeta(2s)}
\frac{\zeta^3(s)}{\zeta(2s)}
=\sum_{n=1}^{\infty}\frac{d(n^2)}{n^s};
\end{align}

\begin{align}
\label{Dirichlet zeta^4(s)/zeta(2s)}
\frac{\zeta^4(s)}{\zeta(2s)}
=\sum_{n=1}^{\infty}\frac{(d(n))^2}{n^s}.
\end{align}

These Dirichlet series involve several arithmetic functions.
The Möbius function ($\mu(n)$) is defined by $\mu(1)=1$; $\mu(n)=(-1)^k$, if $n$ is the product of $k$ distinct primes; and $\mu(n)=0$, if there exists any prime $p$ such that $\dis p^2\mid n$.
The function $\omega(n)$ represents the number of distinct prime factors of $n$.
The Dirichlet divisor function ($d(n)$) expresses the number of divisors of $n$.
For any fixed $k\in\N$, $d_k(n)$ denotes the number of different ways of writing $n$ as the product of $k$ natural factors where expressions with the same factors in different orders are counted as distinct. 
Observe that $d_2(n)=d(n)$.

Finally (see \cite{Ivic}), fixed $t_0>1$, there exists $M\in\R^+$ such that, for all $t\geq t_0$, 
\begin{align}
\label{behaviour of the zeta function on the vertical lines}
\left|\zeta(\sigma\pm it)\right|\leq\begin{cases}
M, & \text{if }\sigma\geq 2;\\
M\ln(t), & \text{if }1\leq\sigma\leq 2;\\
M t^{\frac{1-\sigma}{2}}\ln(t), & \text{if }0\leq\sigma\leq 1;\\
M t^{\frac{1}{2}-\sigma}\ln(t), & \text{if }\sigma\leq 0.
\end{cases}
\end{align}

\newpage
\section{A new family of classes of functions}
The Müntz-type class of functions $\M_\alpha$, where $\alpha>1$, is introduced in \cite{Semyon2}.
Here we generalise this class, defining the following family of classes of functions. 
\begin{definition}
A function $f(x)$, defined for $x\in\R_0^+$, belongs to the generalized Müntz-type class of functions $\M_{\alpha,k}$, where $\alpha>1$ and $k\in\N_0$, if $f\in\mathcal{C}^{(k)}(\R_0^+)$ and $\dis f^{(j)}(x)=\bigO\left(x^{-\alpha-j}\right)$, $x\to\infty$, for all $j=0,1,\cdots,k$.
\end{definition}

This definition is a generalisation of the class $\M_\alpha$ because $\M_{\alpha,2}=\M_{\alpha}$, for any $\alpha>1$. 
Note that, if $k\geq l$ and $\beta\geq\alpha$, then $\M_{\beta,k}\subseteq\M_{\alpha,l}$.

The next theorem shows that the Mellin transform of a function in a $\M_{\alpha,k}$ class is analytic in the strip $\dis -k<\Real s<\alpha$ except some finite (at most $k$) singularity points.
\begin{theorem}
\label{Mellin transform of functions in Malphak}
Let $f\in\M_{\alpha,k}$.
Then, for any $n=0,1,\cdots,k$, $\dis f^{(n)}(x)x^{\sigma+n-1}\in L_1(0,\infty)$, for all $-n<\sigma<\alpha$, and the Mellin transform of $f$, $\dis f^*(s)$, is an analytic function in the strip $0<\Real(s)<\alpha$, which can be analytically continued to the strip $-n<\Real(s)<\alpha$ by
\begin{align}
\label{Mellin transform as integrals over the derivatives of f}
f^*(s)=\frac{(-1)^n}{(s)_n}\int_0^\infty f^{(n)}(x)x^{s+n-1}dx,
\end{align}
where $(s)_n$ is the Pochhammer symbol, defined by $(s)_0=1$ and $(s)_n=s(s+1)\cdots(s+(n-1))$.

Moreover, $f^*(s)$ is analytic in the strip $\dis -k<\Real(s)<\alpha$, \vspace*{0,1 cm} except at the points $s=-n$, with $n=0,1,\cdots,k-1$, where $\dis f^*(s)$ either has a simple pole with residue $\dis\frac{f^{(n)}(0)}{n!}$, if $\dis f^{(n)}(0)\neq 0$, or has a removable singularity, if $\dis f^{(n)}(0)=0$.
\end{theorem}

\textit{Proof.}
Fix $n=0,1,\cdots,k$.
Then, if $-n<\sigma<\alpha$, $\dis f^{(n)}(x)x^{\sigma+n-1}= \bigO\left(x^{\sigma+n-1}\right)$, $x\to 0$ and $\dis f^{(n)}(x)x^{\sigma+n-1}= \bigO\left(x^{\sigma-\alpha-1}\right)$, $x\to\infty$, so $\dis f^{(n)}(x)x^{\sigma+n-1}\in L_1\left(0,\infty\right)$. \vspace*{0,1 cm} 
As a consequence, the integral $\dis\int_{0}^{\infty}f^{(n)}(x)x^{s+n-1}dx$ \vspace*{0,1 cm} defines an analytic function in the strip $-n<\Real(s)<\alpha$.
In particular, if $n=0$, it can be deduced that $\dis f^*(s)$ is analytic in the strip $0<\Real(s)<\alpha$ and its derivatives are obtained differentiating inside the integral (\ref{Mellin transform definition}).

Now we derive (\ref{Mellin transform as integrals over the derivatives of f}) in the strip $0<\Real(s)<\alpha$, for all $n=0,1,\cdots,k$. 
If $n=0$, (\ref{Mellin transform as integrals over the derivatives of f}) coincides with the definition of $\dis f^*(s)$. 
Otherwise, if $n=1,2,\cdots,k$, (\ref{Mellin transform as integrals over the derivatives of f}) can be deduced from the case $n-1$, using integration by parts and eliminating the integrated terms due to the asymptotic behaviour of $f^{(n-1)}(x)$ at the infinity.

Moreover, (\ref{Mellin transform as integrals over the derivatives of f}) gives the analytic continuation of $\dis f^*(s)$ to the strip $-k<\Real(s)<\alpha$, except at the zeros of \vspace*{0,2 cm} $\dis(s)_{k}=s(s+1)\cdots(s+k-1)$: the points $s=-n$, $n=0,1,\cdots,k-1$.
Furthermore, $\dis\lim_{s\to -n}(s+n)f^*(s)
=\lim_{s\to -n}\frac{(-1)^{n+1}(s+n)}{(s)_{n+1}}
 \int_{0}^{\infty}f^{(n+1)}(x)x^{s+n}dx$. \vspace*{0,2 cm}
Besides that, $\dis\frac{(s+n)}{(s)_{n+1}}=\frac{1}{(s)_{n}}$ and $\dis(-n)_n=(-1)^nn!$ \vspace*{0,2 cm} so $\dis\lim_{s\to-n}\frac{(-1)^{n+1}(s+n)}{(s)_{n+1}} =\frac{(-1)^{n+1}}{(-n)_{n}}=-\frac{1}{n!}$
and the integral \vspace*{0,2 cm} $\dis\int_{0}^{\infty}f^{(n+1)}(x)x^{s+n}dx$ defines an analytic function in the strip $-(n+1)<\Real(s)<\alpha$ so $\dis\lim_{s\to-n}\int_{0}^{\infty}f^{(n+1)}(x)x^{s+n}dx =\int_{0}^{\infty}f^{(n+1)}(x)dx=-f^{(n)}(0)$. 
\vspace*{0,1 cm}

Therefore $\dis\lim_{s\to-n}(s+n)f^*(s)=\frac{f^{(n)}(0)}{n!}$ 
\vspace*{0,1 cm}
and, as a result, either $s=-n$ is a simple pole of $\dis f^*(s)$ with residue $\dis\frac{f^{(n)}(0)}{n!}$, if $\dis f^{(n)}(0)\neq 0$, or it is a removable singularity, if $\dis f^{(n)}(0)=0$.\\
\qed

\vfill
The following proposition establishes an upper bound for the Mellin transform of a function in a $\M_{\alpha,k}$ class and, as a result, it gives us sufficient conditions for the absolute convergence of its integral over vertical lines of the complex plane.
\begin{proposition}
\label{Mellin transform of a function in Malphak - upper bound and integrability over vertical lines}
Let $f\in \M_{\alpha,k}$.
\vspace*{0,15 cm}
Then, for any $-k<\sigma<\alpha$, there exists $C(\sigma)\in\R$ such that $\dis|f^*(\sigma+it)|\leq\frac{C(\sigma)}{|t|^k}$, for all $t\in\R\backslash\{0\}$. \vspace*{0,15 cm} 
Moreover, if $k\geq 2$, $\dis f^*(s)\in L_1(\sigma\pm i\infty)$, for any $-k<\sigma<\alpha$ such that $\sigma\neq 0,-1,\cdots,-(k-1)$, and $f(x)$ can be represented by (\ref{Mellin inverse transform}), for any $0<\sigma<\alpha$.
\end{proposition}
\textit{Proof.} 
For any $s\in\C$, $\dis\left|(s)_k\right|\geq|\Imag(s)|^k$ so, if $-k<\Real(s)<\alpha$ and $\Imag(s)\neq 0$, we replace $n=k$ in (\ref{Mellin transform as integrals over the derivatives of f}) to deduce that
\begin{align}
\label{upper bound for f^*(s)}
|f^*(s)|\leq\frac{1}{|(s)_k|}
\int_{0}^{\infty}|f^{(k)}(x)x^{s+k-1}|dx
\leq\frac{1}{|\Imag(s)|^k}
\int_{0}^{\infty}|f^{(k)}(x)|x^{\Real(s)+k-1}dx.
\end{align}

Fix $-k<\sigma<\alpha$. \vspace*{0,15 cm}
Then, by theorem \ref{Mellin transform of functions in Malphak}, $\dis f^{(k)}(x)x^{\sigma+k-1}\in L_1\left(0,\infty\right)$, so we can define \vspace*{0,15 cm}
$\dis C(\sigma)=\int_{0}^{\infty}|f^{(k)}(x)|x^{\sigma+k-1}dx$ and we derive that $\dis|f^*(\sigma+it)|\leq\frac{C(\sigma)}{|t|^k}$, for all $t\in\R\backslash\{0\}$. 
Furthermore, if $\sigma\neq 0,-1,\cdots,-(k-1)$, $\dis f^*(s)$ is continuous on the line $\Real(s)=\sigma$. 
As a result, if $k\geq 2$, $\dis f^*(s)\in L_1(\sigma\pm i\infty)$ and, because $\dis f^*(s)$ is analytic in the strip $0<\Real(s)<\alpha$, $f(x)$ is the Mellin inverse transform of $f^*(s)$ in that strip.\\
\qed

\newpage
\section{Müntz-type formulas in the critical strip}
For functions $f$ with suitable properties (see section 2.11 of \cite{TitchmarshRiemann}), the Müntz formula
\begin{align}
\label{Muntz formula}
\zeta(s)\int_{0}^{\infty}f(y)y^{s-1}dy
=\int_{0}^{\infty}\left(\sum_{n=1}^{\infty}f(nx)
-\frac{1}{x}\int_{0}^{\infty}f(t)dt\right)x^{s-1}dx
\end{align}
is valid in the critical strip $0<\Real(s)<1$ .
Here we derive several identities similar to (\ref{Muntz formula}) in the $\M_{\alpha,k}$ classes.

The following theorem generates, for each Dirichlet series exhibited in the end of our introduction, an equality between an integral over a vertical line in the half-plane $\Real(s)>1$ and a series where appears an arithmetic function (see \cite{Semyon1}).
\begin{theorem}
\label{arithmetic transforms theorem}
Suppose that $f\in\M_{\alpha,k}$, $k\geq 2$, \vspace*{0,1 cm} $\phi(n)$ is an arithmetic function and $\Phi(s)$ is defined in the half-plane $\Real(s)>1$ by the absolutely convergent Dirichlet series $\dis\sum_{n=1}^{\infty}\frac{\phi(n)}{n^s}$. 
Then, for any $1<\sigma<\alpha$, $\dis\Phi(s)f^*(s)\in L_1(\sigma\pm i\infty)$,
\begin{align}
\label{arithmetic transforms - equality between integral and series}
\frac{1}{2\pi i}\int_{\sigma-i\infty}^{\sigma+i\infty}
\Phi(s)f^*(s)x^{-s}ds=\sum_{n=1}^{\infty}\phi(n)f(nx)
\end{align}
and the Mellin transform of $\dis\left(\sum_{n=1}^{\infty}\phi(n)f(nx)\right)$ in the strip $1<\Real(s)<\alpha$ is $\Phi(s)f^*(s)$.
\end{theorem}
\textit{Proof.}
Fix $1<\sigma<\alpha$.
By definition of $\Phi(s)$,
\begin{align*}
\frac{1}{2\pi i}\int_{\sigma-i\infty}^{\sigma+i\infty}
\Phi(s)f^*(s)x^{-s}ds
=\frac{1}{2\pi i}\int_{\sigma-i\infty}^{\sigma+i\infty}
\sum_{n=1}^{\infty}\frac{\phi(n)}{n^s}f^*(s)x^{-s}ds.
\end{align*}

Moreover, the Dirichlet series that defines $\Phi(s)$ is absolutely convergent in the half-plane $\Real(s)>1$ and, by proposition \ref{Mellin transform of a function in Malphak - upper bound and integrability over vertical lines}, $\dis f^*(s)\in L_1(\sigma\pm i\infty)$, so $\dis\Phi(s)f^*(s)\in L_1(\sigma\pm i\infty)$, because
\begin{align*}
\int_{\sigma-i\infty}^{\sigma+i\infty}\left|\Phi(s)f^*(s)ds\right|
\leq\int_{\sigma-i\infty}^{\sigma+i\infty}
\sum_{n=1}^{\infty}\left|\frac{\phi(n)}{n^s}f^*(s)ds\right|
\leq\sum_{n=1}^{\infty}\frac{|\phi(n)|}{n^{\sigma}}
\int_{\sigma-i\infty}^{\sigma+i\infty}\left|f^*(s)ds\right|<\infty.
\end{align*}

Next we change the order of summation and integration to obtain
\begin{align*}
\frac{1}{2\pi i}\int_{\sigma-i\infty}^{\sigma+i\infty}
\Phi(s)f^*(s)x^{-s}ds
=\sum_{n=1}^{\infty}\frac{\phi(n)}{2\pi i} \int_{\sigma-i\infty}^{\sigma+i\infty}f^*(s)(xn)^{-s}ds.
\end{align*}

Furthermore, again by proposition \ref{Mellin transform of a function in Malphak - upper bound and integrability over vertical lines},  $f(x)$ is equal to the inverse Mellin transform of $\dis f^*(s)$ over the vertical line $\Real(s)=\sigma$, so we can deduce \vspace*{0,1 cm} (\ref{arithmetic transforms - equality between integral and series}) from the previous formula.
Finally, using proposition \ref{Mellin transform of a Mellin inverse transform}, we deduce that the Mellin transform of $\dis\left(\sum_{n=1}^{\infty}\phi(n)f(nx)\right)$ exists and is equal to $\Phi(s)f^*(s)$ in the strip $1<\Real(s)<\alpha$.\\
\qed
 
To derive our Müntz-type formulas it will be necessary to move integrals of the type on (\ref{arithmetic transforms - equality between integral and series}), first to the critical strip and later to half-plane $\Real(s)<0$.
To this purpose, we show the following result similar to the residue theorem for integrals over vertical lines of the complex plane.
\begin{theorem}
\label{residue theorem over vertical lines}
Suppose that $a,b,c,d\in\R$, with $c<a<b<d$, $F(s)$ is an analytic function in the strip $c<\Real(s)<d$ except at a point $s=x_0$, with $a<x_0<b$, and there exists $t_0\in\R^+$ and a continuous function $g(t)$ integrable at the infinity such that $\dis\left|F(s)\right|\leq g(\Imag(s))$, for all $s\in\C$ such that $\dis a\leq\Real(s)\leq b$ and $|\Imag(s)|\geq t_0$.

Then $F(s)\in L_1(a\pm i\infty)$, $F(s)\in L_1(b\pm i\infty)$ and, for any $x\in\R^+$,
\begin{align}
\label{change of vertical line of integration with residue}
\res_{s=x_0}\left(F(s)x^{-s}\right)
=\frac{1}{2\pi i}\int_{b-i\infty}^{b+i\infty}F(s)x^{-s}ds
-\frac{1}{2\pi i}\int_{a-i\infty}^{a+i\infty}F(s)x^{-s}ds.
\end{align}
\end{theorem}

\textit{Proof.}
The functions $F(a+it)$ and $F(b+it)$ are continuous and upper bounded at infinity by the integrable function $g(t)$, thus $F(s)\in L_1(a\pm i\infty)$ and $F(s)\in L_1(b\pm i\infty)$.

Fix $x\in\R^+$. For any $T\in\R^+$, we define $\Omega_T$ as the rectangle (positively oriented) whose sides are segments of the vertical lines $\Real(s)=a$ and $\Real(s)=b$ and the horizontal lines $\Imag(s)=T$ and $\Imag(s)=-T$.
Then, using the Cauchy residue theorem, we obtain
\begin{align*}
\res_{s=x_0}\left(F(s)x^{-s}\right)
=\frac{1}{2\pi i}\int_{\Omega_T}F(s)x^{-s}ds
=\frac{1}{2\pi i}\left(\int_{b-iT}^{b+iT}-\int_{a+iT}^{b+iT}
-\int_{a-iT}^{a+iT}+\int_{a-iT}^{b-iT}\right)F(s)x^{-s}ds.
\end{align*}

Besides that, as $g(t)$ is continuous and integrable at the infinity, $\dis\lim_{t\to\infty}g(t)=0$ so
\begin{align*}
\left|\int_{a\pm iT}^{b\pm iT}F(s)x^{-s}ds\right|
\leq\int_{a}^{b}\left|F(u\pm iT)\right|x^{-u}du\leq
\begin{cases}
(b-a)x^{-a}\,g(T), &\text{if }\dis x\geq 1\\
(b-a)x^{-b}\,g(T), &\text{if }\dis x\leq 1
\end{cases}
\xrightarrow{T\to\infty}0.
\end{align*}

Finally we pass to the limit $T\to\infty$ on the formula for $\dis\res_{s=x_0}\left(F(s)x^{-s}\right)$ to deduce (\ref{change of vertical line of integration with residue}).\\
\qed

\subsection{Müntz-type formulas in the critical strip involving $\zeta^k(s)$}
Here we derive a family of Müntz-type formulas where appears $\dis\zeta^k(s)$, $k\in\N$, in the critical strip $0<\Real(s)<1$. 
We observe that the case $k=1$ of these identities is the classical Müntz formula and we determine explicitly the case $k=2$, deducing a Müntz-type formula involving $\dis\zeta^2(s)$.

Fix $k\in\N$ and suppose that $f\in\M_{\alpha,m}$, $\alpha>1$ and $m\in\N_0$. 
Then $\dis\zeta^k(s)f^*(s)$ is analytic in the strip $\dis 0<\Real(s)<\alpha$, except at the point $s=1$ where it has a pole of order at most $k$ (or a removable singularity). 
Moreover, if $m\geq 2$, using (\ref{Dirichlet zeta^k(s)}) and theorem \ref{arithmetic transforms theorem}, we claim that, for all $1<\sigma<\alpha$, $\dis\zeta^k(s)f^*(s)\in L_1(\sigma\pm i\infty)$ and, if $x\in\R^+$,
\begin{align}
\label{arithmetic transform zeta^k(s)}
\frac{1}{2\pi i}\int_{\sigma-i\infty}^{\sigma+i\infty}
\zeta^k(s)f^*(s)x^{-s}ds
=\sum_{n=1}^{\infty}d_k(n)f(nx).
\end{align}

The following proposition establishes sufficient conditions to move the integral in (\ref{arithmetic transform zeta^k(s)}) to the left, using theorem \ref{residue theorem over vertical lines}.
\begin{proposition}
\label{integrability over vertical lines - zeta^k(s)}
Let $k\in\N$ and $f\in\M_{\alpha,m}$, $\dis m\geq 1+\frac{k}{2}$. %\vspace*{0,1 cm}
Then, if we fix $t_0,u_1,u_2\in\R$ such that $t_0>1$ and $\dis\frac{1}{2}-\frac{m-1}{k}<u_1<u_2<\alpha$, there exists a continuous function $g_k(t)$ integrable at the infinity such that $\dis\left|\zeta^k(s)f^*(s)\right|\leq g_k(\Imag(s))$, if $\dis u_1\leq\Real(s)\leq u_2$ and $\dis\left|\Imag(s)\right|\geq t_0$.
\end{proposition}
\vspace*{0,1 cm}

\textit{Proof.}
Fix $t_0,u_1,u_2\in\R$ such that $t_0>1$ and  $\dis\frac{1}{2}-\frac{m-1}{k}<u_1<u_2<\alpha$. \vspace*{0,1 cm} 

For any $u_1\leq u\leq u_2$,
\begin{align*}
\int_{0}^{\infty}\left|f^{(m)}(x)\right|x^{u+m-1}dx
\leq\int_{0}^{1}\left|f^{(m)}(x)\right|x^{u_1+m-1}dx
+\int_{1}^{\infty}\left|f^{(m)}(x)\right|x^{u_2+m-1}dx=:C,
\end{align*}
so, remembering (\ref{upper bound for f^*(s)}), we obtain, for all $t>0$,
\begin{align*}
\left|f^*(u\pm it)\right|
\leq\frac{1}{t^m}\int_{0}^{\infty}|f^{(n)}(x)|x^{u+m-1}dx
\leq\frac{C}{t^m}.
\end{align*}

By (\ref{behaviour of the zeta function on the vertical lines}), there exists $M\in\R^+$ such that, for all $t\geq t_0$:
$\dis\left|\zeta(u\pm it)\right|\leq M\ln(t)$, 
if $u\geq 1$;
$\dis\left|\zeta(u\pm it)\right|\leq M t^{\frac{1-u}{2}}\ln(t)$, if $0\leq u\leq 1$; 
and $\dis\left|\zeta(u\pm it)\right|\leq Mt^{\frac{1}{2}-u}\ln(t)$,
if $u\leq 0$.

Therefore, for any $t\geq t_0$ and $u_1\leq u\leq u_2$,
\vspace*{0,15 cm} $\dis\left|\zeta^k(u\pm it)f^*(u\pm it)\right|\leq g_k(t)=M^kC(\ln(t))^k\,t^{p(u_1)}$,
where $\dis p(u_1)=-m$, if $\dis 1\leq u_1<\alpha$;
$\dis p(u_1)=-m+\frac{k}{2}\left(1-u_1\right)$, 
if $\dis 0\leq u_1\leq 1$; 
and $\dis p(u_1)=-m+k\left(\frac{1}{2}-u_1\right)$, 
if $\dis\frac{1}{2}-\frac{m-1}{k}<u_1\leq 0$. \vspace*{0,1 cm}
Finally, $\dis p(u_1)<-1$, for any of the possible values for $u_1$, so $\dis g_k(t)$ is a (continuous) function integrable at the infinity.\\
\qed

\newpage
Suppose that $f\in\M_{\alpha,m}$, $\dis m\geq 1+\frac{k}{2}$.
\vspace*{0,1 cm}
Fix $\dis 0<c_0<1$ and $1<\sigma<\alpha$.
Then theorem \ref{residue theorem over vertical lines} and proposition \ref{integrability over vertical lines - zeta^k(s)} can be used to deduce that $\dis\zeta^k(s)f^*(s)\in L_1(c_0\pm i\infty)$ and
\begin{align*}
\frac{1}{2\pi i}\int_{c_0-i\infty}^{c_0+i\infty}
\zeta^k(s)f^*(s)x^{-s}ds
=\frac{1}{2\pi i}\int_{\sigma-i\infty}^{\sigma+i\infty}
\zeta^k(s)f^*(s)x^{-s}ds
-\res_{s=1}\left(\zeta^k(s)f^*(s)x^{-s}\right).
\end{align*}

Besides that, 
\begin{align*}
\res_{s=1}\left(\zeta^k(s)f^*(s)x^{-s}\right)
=\int_{0}^{\infty}f(xy)P_{k-1}(\ln(y))dy,
\end{align*}
where $\dis P_{k-1}(x)$ is a certain monic polynomial of degree $k-1$ (see \cite{Semyon2}). \vspace*{0,1 cm}

Therefore, remembering (\ref{arithmetic transform zeta^k(s)}), we obtain
\begin{align}
\label{integral in the critical strip - zeta^k(s)}
\frac{1}{2\pi i}\int_{c_0-i\infty}^{c_0+i\infty}
\zeta^k(s)f^*(s)x^{-s}ds
=\sum_{n=1}^{\infty}d_k(n)f(nx)
-\int_{0}^{\infty}f(xy)P_{k-1}(\ln(y))dy.
\end{align}

Finally, applying proposition \ref{Mellin transform of a Mellin inverse transform} to (\ref{integral in the critical strip - zeta^k(s)}), we derive the following theorem.

\begin{theorem}
\label{Muntz-type formula - zeta^k(s) - theorem}
Let $k\in\N$ and $f\in\M_{\alpha,m}$, $\dis m\geq 1+\frac{k}{2}$.
Then the Müntz-type formula
\begin{align}
\label{Muntz-type formula - zeta^k(s)}
\zeta^k(s)f^*(s)
=\int_{0}^{\infty}\left(\sum_{n=1}^{\infty}d_k(n)f(nx) -\int_{0}^{\infty}f(xy)P_{k-1}(\ln(y))dy\right)x^{s-1}dx
\end{align}
is valid in the critical strip $\dis 0<\Real(s)<1$. 
\end{theorem}

If $k=1$, $\dis P_0(x)=1$ so
\begin{align}
\label{residue zeta(s)f^*(s)x^-s, s=1}
\res_{s=1}\left(\zeta(s)f^*(s)x^{-s}\right)
=\int_{0}^{\infty}f(xy)dy
=\frac{1}{x}\int_{0}^{\infty}f(t)dt
=\frac{f^*(1)}{x}.
\end{align}

Therefore, because $d_1(n)=1$, for all $n\in\N$, we deduce from the previous theorem that the Müntz formula (\ref{Muntz formula}) is valid in the critical strip, for any function $f\in\M_{\alpha,m}$, $m\geq 2$.

If $k=2$, $\dis P_1(x)=x+2\gamma$ (where $\gamma$ is the Euler-Mascheroni constant \cite{Havil}) so
\begin{align*}
\res_{s=1}\left(\zeta^2(s)f^*(s)x^{-s}\right)
=\int_{0}^{\infty}f(xy)P_1(\ln(y))dy
=2\gamma\int_{0}^{\infty}f(xy)dy+\int_{0}^{\infty}f(xy)\ln(y)dy.
\end{align*}

Besides that, making a change of variable $t=xy$,
\begin{align*}
\int_{0}^{\infty}f(xy)\ln(y)dy
=\frac{1}{x}\int_{0}^{\infty}f(t)\ln(t)dt
-\frac{\ln(x)}{x}\int_{0}^{\infty}f(t)dt
=\frac{(f^*)'(1)}{x}-f^*(1)\frac{\ln(x)}{x}.
\end{align*}

As a result,
\begin{align}
\label{residue zeta^2(s)f^*(s)x^-s, s=1}
\res_{s=1}\left(\zeta^2(s)f^*(s)x^{-s}\right)
=\int_{0}^{\infty}f(xy)P_1(\ln(y))dy
=\frac{1}{x}\,\Big(\big((f^*)'(1)+2\gamma f^*(1)\big)-f^*(1)\ln(x)\Big)
\end{align}
and, as $d_2(n)=d(n)$, for all $n\in\N$, we derive the Müntz-type formula
\begin{align}
\label{Muntz-type formula - zeta^2(s)}
\zeta^2(s)f^*(s)
=\int_{0}^{\infty}\left(\sum_{n=1}^{\infty}d(n)f(nx)+\frac{1}{x}\,
\Big(f^*(1)\ln(x)-\big((f^*)'(1)+2\gamma f^*(1)\big)\Big)\right)x^{s-1}dx
\end{align}
valid in the critical strip $\dis 0<\Real(s)<1$, for any function $f\in\M_{\alpha,m}$, $m\geq 2$.

%\newpage
\subsection{Müntz-type formulas in the critical strip involving $\frac{\zeta^k(s)}{\zeta(2s)}$}
In this section we derive Müntz-type formulas where appears $\dis\frac{\zeta^k(s)}{\zeta(2s)}$, $k=1,2,3,4$. \vspace*{0,2 cm}

Let $k\in\N$ and $f\in\M_{\alpha,m}$. \vspace*{0,1 cm}
Then, because $\dis\frac{1}{\zeta(2s)}$ is analytic in the half-plane $\dis\Real(s)>\frac{1}{2}$, $\dis\frac{\zeta^k(s)}{\zeta(2s)}f^*(s)$  \vspace*{0,1 cm} is an analytic function in the strip $\dis\frac{1}{2}<\Real(s)<\alpha$ except at the point $s=1$, where it has a pole of order at most $k$ (or a removable singularity).

Moreover, if $m\geq 2$, remembering formulas (\ref{Dirichlet zeta(s)/zeta(2s)})-(\ref{Dirichlet zeta^4(s)/zeta(2s)}) and theorem \ref{arithmetic transforms theorem}, \vspace*{0,1 cm} 
we deduce that, for all $1<\sigma<\alpha$, $\dis\frac{\zeta^k(s)}{\zeta(2s)}f^*(s)\in L_1(\sigma\pm i\infty)$ and, for any $x\in\R^+$,
\begin{align}
\label{arithmetic transform zeta(s)/zeta(2s)}
\frac{1}{2\pi i}\int_{\sigma-i\infty}^{\sigma+i\infty}
\frac{\zeta(s)}{\zeta(2s)}f^*(s)x^{-s}ds
=\sum_{n=1}^{\infty}|\mu(n)|f(nx);
\end{align}

\begin{align}
\label{arithmetic transform zeta^2(s)/zeta(2s)}
\frac{1}{2\pi i}\int_{\sigma-i\infty}^{\sigma+i\infty}
\frac{\zeta^2(s)}{\zeta(2s)}f^*(s)x^{-s}ds
=\sum_{n=1}^{\infty}2^{\omega(n)}f(nx);
\end{align}

\begin{align}
\label{arithmetic transform zeta^3(s)/zeta(2s)}
\frac{1}{2\pi i}\int_{\sigma-i\infty}^{\sigma+i\infty}
\frac{\zeta^3(s)}{\zeta(2s)}f^*(s)x^{-s}ds
=\sum_{n=1}^{\infty}d(n^2)f(nx);
\end{align}

\begin{align}
\label{arithmetic transform zeta^4(s)/zeta(2s)}
\frac{1}{2\pi i}\int_{\sigma-i\infty}^{\sigma+i\infty}
\frac{\zeta^4(s)}{\zeta(2s)}f^*(s)x^{-s}ds
=\sum_{n=1}^{\infty}(d(n))^2f(nx).
\end{align}

Let $u_1\in\R$ and $s\in\C$ such that $\dis u:=\Real(s)\geq u_1>\frac{1}{2}$. 
Then, using (\ref{zeta function definition}) and (\ref{Dirichlet 1/zeta}), 
\begin{align*}
\left|\frac{1}{\zeta(2s)}\right|
\leq\sum_{n=1}^{\infty}\left|\frac{\mu(n)}{n^{2s}}\right|
\leq\sum_{n=1}^{\infty}\frac{1}{n^{2u}}
=\zeta(2u)\leq\zeta(2u_1).
\end{align*}

Therefore $\dis\left|\frac{\zeta^k(s)}{\zeta(2s)}f^*(s)\right|
\leq\zeta(2u_1)\left|\zeta^k(s)f^*(s)\right|$ \vspace*{0,1 cm}
and we may derive the following result from proposition \ref{integrability over vertical lines - zeta^k(s)}, defining $\dis\tilde{g}_k(t):=\zeta(2u_1)g_k(t)$, where $g_k(t)$ is the function obtained in proposition \ref{integrability over vertical lines - zeta^k(s)}.
\begin{proposition}
\label{integrability over vertical lines - zeta^k(s)/zeta(2s)}
Let $k\in\N$ and $f\in\M_{\alpha,m}$, $\dis m\geq 1+\frac{k}{2}$. \vspace*{0,15 cm}
Then, if we fix $t_0,u_1,u_2\in\R$ such that $t_0>1$ and $\dis\frac{1}{2}<u_1<u_2<\alpha$, \vspace*{0,15 cm} there exists a continuous function $\dis\tilde{g}_k(t)$ integrable at the infinity such that $\dis\left|\frac{\zeta^k(s)}{\zeta(2s)}f^*(s)\right|\leq \dis\tilde{g}_k(\Imag(s))$, if $\dis u_1\leq\Real(s)\leq u_2$ and $\left|\Imag(s)\right|\geq t_0$.
\end{proposition}

%\newpage
Suppose that $f\in\M_{\alpha,m}$, $\dis m\geq 1+\frac{k}{2}$.
Fix $\dis\frac{1}{2}<c_0<1$ and $1<\sigma<\alpha$. \vspace*{0,1 cm}
Then theorem \ref{residue theorem over vertical lines} and proposition \ref{integrability over vertical lines - zeta^k(s)/zeta(2s)} can be used to deduce that $\dis\frac{\zeta^k(s)}{\zeta(2s)}f^*(s)\in L_1(c_0\pm i\infty)$ and
\begin{align}
\label{integral in the critical strip - zeta^k(s)/zeta(2s)}
\frac{1}{2\pi i}\int_{c_0-i\infty}^{c_0+i\infty}
\frac{\zeta^k(s)}{\zeta(2s)}f^*(s)x^{-s}ds
=\frac{1}{2\pi i}\int_{\sigma-i\infty}^{\sigma+i\infty}
\frac{\zeta^k(s)}{\zeta(2s)}f^*(s)x^{-s}ds
-\res_{s=1}\left(\frac{\zeta^k(s)}{\zeta(2s)}f^*(s)x^{-s}\right).
\end{align}

Next we calculate this residue for $k=1$ and $k=2$. 
\vspace*{0,1 cm} 

If $k=1$, $\dis\frac{\zeta(s)}{\zeta(2s)}f^*(s)x^{-s}$ has a simple pole (or a removable singularity) at the point $s=1$ and then, remembering (\ref{residue zeta(s)f^*(s)x^-s, s=1}) and $\dis\zeta(2)=\frac{\pi^2}{6}$,
\begin{align}
\label{residue zeta(s)/zeta(2s)f^*(s)x^-s, s=1}
\res_{s=1}\left(\frac{\zeta(s)}{\zeta(2s)}f^*(s)x^{-s}\right)
=\frac{1}{\zeta(2)}\res_{s=1}\left(\zeta(s)f^*(s)x^{-s}\right)
=\frac{6f^*(1)}{\pi^2 x}.
\end{align}

If $k=2$, $\dis\frac{\zeta^2(s)}{\zeta(2s)}f^*(s)x^{-s}$ has (at most) a double pole at the point $s=1$, so
\begin{align*}
&\res_{s=1}\left(\frac{\zeta^2(s)}{\zeta(2s)}f^*(s)x^{-s}\right)
=\lim_{s\to 1}\frac{d}{ds}
\left((s-1)^2\frac{\zeta^2(s)}{\zeta(2s)}f^*(s)x^{-s}\right)\\
&=\frac{1}{\zeta(2)}\res_{s=1}\left(\zeta^2(s)f^*(s)x^{-s}\right) +\frac{1}{x}\,f^*(1)\left(\lim_{s\to 1}\big((s-1)\zeta(s)\big)\right)^2
\left.\frac{d}{ds}\left(\frac{1}{\zeta(2s)}\right)\right|_{s=1}.
\end{align*}

Besides that $\dis\zeta'(2)=\frac{\pi^2}{6} \left(\gamma+\ln\left(\frac{2\pi}{A^{12}}\right)\right)$, where $A$ is the Glaisher-Kinkelin constant, so
\begin{align*}
\left.\frac{d}{ds}\left(\frac{1}{\zeta(2s)}\right)\right|_{s=1}
=-\frac{2\zeta'(2)}{\zeta^2(2)}
=\frac{12}{\pi^2}\left(\ln\left(\frac{A^{12}}{2\pi}\right)-\gamma\right).
\end{align*}

Therefore, remembering (\ref{residue zeta^2(s)f^*(s)x^-s, s=1}), we obtain
\begin{align}
\label{residue zeta^2(s)/zeta(2s)f^*(s)x^-s, s=1}
\res_{s=1}\left(\frac{\zeta^2(s)}{\zeta(2s)}f^*(s)x^{-s}\right)
=\frac{6}{\pi^2x}
\left((f^*)'(1)+f^*(1)\ln\left(\frac{A^{24}}{4\pi^2 x}\right)\right).
\end{align}

Finally, the following theorem can be derived from (\ref{integral in the critical strip - zeta^k(s)/zeta(2s)}), (\ref{arithmetic transform zeta(s)/zeta(2s)})-(\ref{arithmetic transform zeta^2(s)/zeta(2s)}), (\ref{residue zeta(s)/zeta(2s)f^*(s)x^-s, s=1})-(\ref{residue zeta^2(s)/zeta(2s)f^*(s)x^-s, s=1}) and proposition \ref{Mellin transform of a Mellin inverse transform}.
\begin{theorem}{}
\label{Muntz-type formula - zeta^2(s)/zeta(2s)}
Let $f\in\M_{\alpha,m}$, $m\geq 2$. \vspace*{0,1 cm} 
Then the Müntz-type formulas
\begin{align}
\frac{\zeta(s)}{\zeta(2s)}f^*(s)
=\int_{0}^{\infty}\left(\sum_{n=1}^{\infty}|\mu(n)|f(nx)
-\frac{6f^*(1)}{\pi^2 x}\right)x^{s-1}dx,
\end{align}
and 
\begin{align}
\frac{\zeta^2(s)}{\zeta(2s)}f^*(s)
=\int_{0}^{\infty}\left(\sum_{n=1}^{\infty}2^{\omega(n)}f(nx)
+\frac{6}{\pi^2x}\left(f^*(1)\ln\left(\frac{4\pi^2x}{A^{24}}\right) -(f^*)'(1)\right)\right)x^{s-1}dx
\end{align}
are valid in the strip $\dis\frac{1}{2}<\Real(s)<1$.
\end{theorem}

Analogously, using (\ref{arithmetic transform zeta^3(s)/zeta(2s)}), (\ref{arithmetic transform zeta^4(s)/zeta(2s)}) and (\ref{integral in the critical strip - zeta^k(s)/zeta(2s)}), one may deduce the following theorem (but we will not calculate the residues appearing here).
\begin{theorem}
\label{more Muntz-type formulas}
Let $f\in\M_{\alpha,m}$, $m\geq 3$. \vspace*{0,1 cm} 
Then the Müntz-type formulas
\begin{align}
\label{Muntz-type formula - zeta^3(s)/zeta(2s)}
\frac{\zeta^3(s)}{\zeta(2s)}f^*(s)
=\int_{0}^{\infty}\left(\sum_{n=1}^{\infty}d(n^2)f(nx)
-\res_{s=1}\left(\frac{\zeta^3(s)}{\zeta(2s)}f^*(s)x^{-s}\right)
\right)x^{s-1}dx
\end{align}
and
\begin{align}
\label{Muntz-type formula - zeta^4(s)/zeta(2s)}
\frac{\zeta^4(s)}{\zeta(2s)}f^*(s)
=\int_{0}^{\infty}\left(\sum_{n=1}^{\infty}d^2(n)f(nx)
-\res_{s=1}\left(\frac{\zeta^4(s)}{\zeta(2s)}f^*(s)x^{-s}\right)
\right)x^{s-1}dx
\end{align}
are valid in the strip $\dis\frac{1}{2}<\Real(s)<1$.
\end{theorem}

This section ends with a remark about how the Riemann hypothesis may affect the strip of validity of the formulas exhibited above. 

\begin{remark}
If the Riemann hypothesis holds true, \vspace*{0,1 cm} the formulas given by theorems \ref{Muntz-type formula - zeta^2(s)/zeta(2s)} and \ref{more Muntz-type formulas} are not only valid in the strip $\dis\frac{1}{2}<\Real(s)<1$ but in the entire strip $\dis\frac{1}{4}<\Real(s)<1$. 
\end{remark}

\newpage
\section{Müntz-type formulas in the half-plane $\Real(s)<0$}
In the previous section we moved some integrals of the type on (\ref{arithmetic transforms - equality between integral and series}) to the critical strip $0<\Real(s)<1$. 
Now we move the integrals of $\dis\zeta^k(s)f^*(s)x^{-s}$ ($k\in\N$) to the half-plane $\Real(s)<0$ and we derive some Müntz-type formulas in that half-plane.

Suppose that $f\in\M_{\alpha,m}$, $\alpha>1$ and $m\in\N$.
By theorem \ref{Mellin transform of functions in Malphak}, $\dis f^*(s)$ is analytic in the strip $-m<\Real(s)<\alpha$ except at the points $s=-j$, $j=0,1,\cdots,m-1$, where $\dis f^*(s)$ either has a simple pole, if $f^{(j)}(0)\neq 0$, or a removable singularity, if $f^{(j)}(0)=0$.
Besides that, $\dis\zeta(s)$ is analytic in the entire complex plane except at the point $s=1$ and $\zeta(-2n)=0$, for all $j\in\N$, so these zeros cancel the eventual simple poles of $\dis f^*(s)$ at the points $s=-2n$, $1<2n<m$.
Therefore, for any fixed $k\in\N$, $\dis\zeta^k(s)f^*(s)x^{-s}$ is an analytic function in the strip $-m<\Real(s)<1$, except at the points $s=-j$, with $j=0$ or $j=2n-1<m$ ($n\in\N$), where $\dis\zeta^k(s)f^*(s)x^{-s}$ either has a simple pole, if $\dis f^{(j)}(0)\neq 0$, or a removable singularity, if $\dis f^{(j)}(0)=0$.
Moreover, we can use theorem \ref{Mellin transform of functions in Malphak} to calculate the residues of $\dis\zeta^k(s)f^*(s)x^{-s}$ at these possible singularities, because
\begin{align}
\label{residues of zeta^k(s)f^*(s)x^-s, s<0}
\res_{s=-j}\left(\zeta^k(s)f^*(s)x^{-s}\right)
=\zeta^k(-j)\,x^{j}\res_{s=-j}f^*(s)
=\frac{\zeta^k(-j)}{j!}\,f^{(j)}(0)\,x^{j},
\end{align}
for all $j\in\N_0$ such that $j<m$.
In particular, if $j=0$, then, because $\dis\zeta(0)=-\frac{1}{2}$,
\begin{align}
\label{residue of zeta^k(s)f^*(s)x^-s at s=0}
\res_{s=0}\left(\zeta^k(s)f^*(s)x^{-s}\right)
=\frac{(-1)^k}{2^k}\,f(0).
\end{align}

Now we move the integral of $\zeta^k(s)f^*(s)x^{-s}$ to the strip $-1<\Real(s)<0$ to obtain a Müntz-type formula in that strip. 

Let $k\in\N$ and $f\in\M_{\alpha,m}$, $\dis m\geq 1+\frac{3k}{2}$.
Fix $-1<\sigma_0<0$ and $0<c_0<1$. \vspace*{0,1 cm}
Using proposition \ref{integrability over vertical lines - zeta^k(s)} and  theorem \ref{residue theorem over vertical lines}, we claim that $\dis\zeta^k(s)f^*(s)\in L_1(\sigma_0\pm i\infty)$ and
\begin{align}
\frac{1}{2\pi i}\int_{\sigma_0-i\infty}^{\sigma_0+i\infty}
\zeta^k(s)f^*(s)x^{-s}ds
=\frac{1}{2\pi i}\int_{c_0-i\infty}^{c_0+i\infty}
\zeta^k(s)f^*(s)x^{-s}ds
-\res_{s=0}\left(\zeta^k(s)f^*(s)x^{-s}\right).
\end{align}

Then, remembering (\ref{integral in the critical strip - zeta^k(s)}) and (\ref{residue of zeta^k(s)f^*(s)x^-s at s=0}), we deduce that
\begin{align}
\label{integral of zeta^k(s)f*(s)x^-s in the strip -1<Re s<0}
\frac{1}{2\pi i}\int_{\sigma_0-i\infty}^{\sigma_0+i\infty}
\zeta^k(s)f^*(s)x^{-s}ds
=\sum_{n=1}^{\infty}d_k(n)f(nx)
-\int_{0}^{\infty}f(xy)P_{k-1}(\ln(y))dy
+\frac{(-1)^{k+1}}{2^k}\,f(0).
\end{align}

\newpage
Therefore, the following theorem can be obtained applying proposition \ref{Mellin transform of a Mellin inverse transform} to (\ref{integral of zeta^k(s)f*(s)x^-s in the strip -1<Re s<0}).
\begin{theorem}
\label{Muntz-type formula - -1<Re s<0}
Let $k\in\N$ and $f\in\M_{\alpha,m}$, $\dis m\geq 1+\frac{3k}{2}$. Then the Müntz-type formula
\begin{align}
\zeta^k(s)f^*(s)
=\int_{0}^{\infty}\left(\sum_{n=1}^{\infty}d_k(n)f(nx)
-\int_{0}^{\infty}f(xy)P_{k-1}(\ln(y))dy
+\frac{(-1)^{k+1}}{2^k}\,f(0)\right)x^{s-1}dx
\end{align}
is valid in the strip $-1<\Real(s)<0$.
\end{theorem}

Replacing $k=1$ and $k=2$ on the previous theorem and remembering (\ref{residue zeta(s)f^*(s)x^-s, s=1}) and (\ref{residue zeta^2(s)f^*(s)x^-s, s=1}), we derive the Müntz-type formulas
\begin{align}
\label{Muntz-type formula - -1<Re s<0 - k=1}
\zeta(s)f^*(s)=\int_{0}^{\infty}
\left(\sum_{n=1}^{\infty}f(nx)-\frac{f^*(1)}{x}
+\frac{f(0)}{2}\right)x^{s-1}dx
\end{align}
and 
\begin{align}
\label{Muntz-type formula - -1<Re s<0 - k=2}
\zeta^2(s)f^*(s)=\int_{0}^{\infty}\left(\sum_{n=1}^{\infty}d(n)f(nx) +\frac{1}{x}\,\Big(f^*(1)\ln(x)-\big((f^*)'(1)+2\gamma f^*(1)\big)\Big) -\frac{f(0)}{4}\right)x^{s-1}dx
\end{align}
valid in the strip $-1<\Real(s)<0$, for any function $f\in\M_{\alpha,m}$, where $m\geq 3$ and $m\geq 4$, respectively.

Next we move the integral of $\zeta^k(s)f^*(s)x^{-s}$ to the half-plane $\Real(s)<-1$ to deduce Müntz-type formulas in strips of that half-plane.
\vspace*{0,1 cm}

Let $k,m\in\N$ and $f\in\M_{\alpha,l}$, $\dis l\geq 1+k\left(2m+\frac{3}{2}\right)$. \vspace*{0,2 cm}
Fix $\sigma_n$, $n=0,1,\cdots,m$, such that $-1<\sigma_0<0$ and $-2n-1<\sigma_n<-2n+1$, for all $n=1,2,\cdots,m$.
Using proposition \ref{integrability over vertical lines - zeta^k(s)} and theorem \ref{residue theorem over vertical lines}, we claim that, for any $n=1,2,\cdots,m$, $\dis\zeta^k(s)f^*(s)\in L_1(\sigma_n\pm i\infty)$ and
\begin{align*}
\frac{1}{2\pi i}\int_{\sigma_{n}-i\infty}^{\sigma_{n}+i\infty}
\zeta^k(s)f^*(s)x^{-s}ds
=\frac{1}{2\pi i}
\int_{\sigma_{n-1}-i\infty}^{\sigma_{n-1}+i\infty}
\zeta^k(s)f^*(s)x^{-s}ds
-\res_{s=1-2n}\left(\zeta^k(s)f^*(s)x^{-s}\right).
\end{align*}

Then, remembering formulas (\ref{residues of zeta^k(s)f^*(s)x^-s, s<0}) and (\ref{integral of zeta^k(s)f*(s)x^-s in the strip -1<Re s<0}), we may deduce that
\begin{align}
\label{integral of zeta^k(s)f*(s)x^-s, Re s<-1}
\frac{1}{2\pi i}\int_{\sigma_m-i\infty}^{\sigma_m+i\infty}
\zeta^k(s)f^*(s)x^{-s}ds=(P_{k,m}f)(x),
\end{align}
where $\dis (P_{k,m}f)(x)$ ($x\in\R^+$) is defined as
\begin{align*}
\sum_{n=1}^{\infty}d_k(n)f(nx)
-\int_{0}^{\infty}f(xy)P_{k-1}(\ln(y))dy
+\frac{(-1)^{k+1}f(0)}{2^k}
-\sum_{n=1}^{m}\frac{\zeta^k(1-2n)}{(2n-1)!}\,f^{(2n-1)}(0)\,x^{2n-1}.
\end{align*} 

\newpage
Finally, the following theorem can be obtained applying proposition \ref{Mellin transform of a Mellin inverse transform} to (\ref{integral of zeta^k(s)f*(s)x^-s, Re s<-1}).
\begin{theorem}
\label{Muntz-type formula - Re s<-1}
Let $k,m\in\N$ and $f\in\M_{\alpha,l}$, $\dis l\geq 1+k\left(2m+\frac{3}{2}\right)$.
Then, if we define $\dis(P_{k,m}f)(x)$ as above, the Müntz-type formula
\begin{align}
\zeta^k(s)f^*(s)=\int_{0}^{\infty}(P_{k,m}f)(x)x^{s-1}dx
\end{align}
is valid in the strip $-2m-1<\Real(s)<-2m+1$.
\end{theorem}

Fixing $m\in\N$, replacing $k=1$ and $k=2$ on the previous theorem and remembering (\ref{residue zeta(s)f^*(s)x^-s, s=1}) and (\ref{residue zeta^2(s)f^*(s)x^-s, s=1}), we derive the Müntz-type formulas
\begin{align}
\label{Muntz-type formula - Re s<-1 - k=1}
\zeta(s)f^*(s)=\int_{0}^{\infty}
\left(\sum_{n=1}^{\infty}f(nx)-\frac{f^*(1)}{x}+\frac{f(0)}{2}
-\sum_{n=1}^{m}\frac{\zeta(1-2n)}{(2n-1)!}\,f^{(2n-1)}(0)\,x^{2n-1}
\right)x^{s-1}dx
\end{align}
and
\begin{align}
\label{Muntz-type formula - Re s<-1 - k=2}
\zeta^2(s)f^*(s)=\int_{0}^{\infty}(P_{2,m}f)(x)x^{s-1}dx,
\end{align}
where $\dis (P_{2,m}f)(x)$ ($x\in\R^+$) is defined as equal to
\begin{align*} 
\sum_{n=1}^{\infty}d(n)f(nx)
+\frac{1}{x}\,\Big(f^*(1)\ln(x)-\big((f^*)'(1)+2\gamma f^*(1)\big)\Big) -\frac{f(0)}{4}
-\sum_{n=1}^{m}\frac{\zeta^2(1-2n)}{(2n-1)!}\,f^{(2n-1)}(0)\,x^{2n-1},
\end{align*}
valid in the strip $-2m-1<\Real(s)<-2m+1$, for any function $f\in\M_{\alpha,l}$, where $l\geq 2m+3$ and $l\geq 4m+4$, respectively.
\vspace*{0,2 cm}

This section ends with a remark about the relation between the Müntz-type formulas in the strip $-1<\Real(s)<0$ and some classical summation formulas.
\begin{remark}
The classical Poisson and Voronoi summation formulas (see \cite{TitchmarshFourier} and \cite{Ivic}, respectively) can be derived in the $\M_{\alpha,m}$ classes, replacing $k=1$ and $k=2$, respectively, on (\ref{integral of zeta^k(s)f*(s)x^-s in the strip -1<Re s<0}) \vspace*{0,15 cm} and using the functional equation of the Riemann zeta function \cite{TitchmarshRiemann} to calculate $\dis\frac{1}{2\pi i}\int_{\sigma_0-i\infty}^{\sigma_0+i\infty} \zeta^k(s)f^*(s)x^{-s}ds$, $\dis-1<\sigma_0<0$.
\end{remark}

\newpage
\section{Identities involving the gamma and zeta functions}
Perhaps the most important example of a function belonging to the generalized Müntz-type classes of functions is $f(x)=e^{-x}\in\M_{\alpha,k}$, for all $\alpha>1$ and $k\in\N_0$.
Its Mellin transform is the gamma function $\Gamma(s)$.
In this section we replace $f(x)=e^{-x}$ and $\dis f^*(s)=\Gamma(s)$ in the previously derived formulas to obtain integral representations in vertical strips of the complex plane for products of the gamma and zeta functions.

Observe that, if $f(x)=e^{-x}$ and $\dis f^*(s)=\Gamma(s)$, then $\dis f^*(1)=\Gamma(1)=1$ and
\begin{align*}
\sum_{n=1}^{\infty}f(nx)
=\sum_{n=1}^{\infty}e^{-nx}=\frac{1}{e^x-1}.
\end{align*}

Therefore we can derive from the Müntz formula the following integral representation for $\dis\zeta(s)\Gamma(s)$ in the critical strip $\dis 0<\Real(s)<1$, which may be found in section 2.7 of \cite{TitchmarshRiemann},
\begin{align}
\zeta(s)\Gamma(s)=\int_{0}^{\infty}
\left(\frac{1}{e^x-1}-\frac{1}{x}\right)x^{s-1}dx.
\end{align}

Similarly, we may deduce more formulas relating the gamma and zeta functions.
For instance, making $f(x)=e^{-x}$ and $\dis f^*(s)=\Gamma(s)$, \vspace*{0,1 cm} we have $\dis f^*(1)=1$ and $\dis(f^*)'(1)=\Gamma'(1)=-\gamma$, so $\dis f^*(1)\,\frac{\ln(x)}{x}-\big(2\gamma f^*(1)+(f^*)'(1)\big)\frac{1}{x}=\frac{\ln(x)}{x}-\frac{\gamma}{x}$. 
As a result, we obtain from (\ref{Muntz-type formula - zeta^2(s)}) the integral representation
\begin{align}
\zeta^2(s)\Gamma(s)
=\int_{0}^{\infty}x^{s-1}\left(\sum_{n=1}^{\infty}d(n)e^{-nx}
+\frac{\ln(x)}{x}-\frac{\gamma}{x}\right)dx
\end{align}
valid in the crtical strip $\dis 0<\Real(s)<1$.

Analogously, we derive from theorem \ref{Muntz-type formula - zeta^2(s)/zeta(2s)} the integral representations
\begin{align}
\label{formula for zeta(s)/zeta(2s)*Gamma(s)}
\frac{\zeta(s)}{\zeta(2s)}\Gamma(s)
=\int_{0}^{\infty}x^{s-1}
\left(\sum_{n=1}^{\infty}|\mu(n)|e^{-nx}
-\frac{6}{\pi^2 x}\right)dx
\end{align}
and
\begin{align}
\label{formula for zeta^2(s)/zeta(2s)*Gamma(s)}
\frac{\zeta^2(s)}{\zeta(2s)}\Gamma(s)
=\int_{0}^{\infty}x^{s-1}
\left(\sum_{n=1}^{\infty}2^{\omega(n)}e^{-nx}
+\frac{6}{\pi^2x}\left(\ln\left(\frac{4\pi^2x}{A^{24}}\right)
+\gamma\right)\right)dx
\end{align}
valid in the strip $\dis\frac{1}{2}<\Real(s)<1$
(or in the strip $\dis\frac{1}{4}<\Real(s)<1$, if the Riemann hypothesis holds true).

\newpage
Moreover, if $f(x)=e^{-x}$ then $f(0)=1$, so we deduce from (\ref{Muntz-type formula - -1<Re s<0 - k=1}) and (\ref{Muntz-type formula - -1<Re s<0 - k=2}) the integral representations
\begin{align}
\zeta(s)\Gamma(s)=\int_{0}^{\infty}
\left(\frac{1}{e^x-1}-\frac{1}{x}+\frac{1}{2}\right)x^{s-1}dx
\end{align}
and
\begin{align}
\zeta^2(s)\Gamma(s) =\int_{0}^{\infty}\left(\sum_{n=1}^{\infty}d(n)e^{-nx} +\frac{\ln(x)}{x}-\frac{\gamma}{x}-\frac{1}{4}\right)x^{s-1}dx
\end{align}
valid in the strip $\dis -1<\Real(s)<0$.

Finally, if $f(x)=e^{-x}$ then, for any $j\in\N_0$, $f^{(j)}(x)=(-1)^je^{-x}$ so $f^{(j)}(0)=(-1)^j$.
In particular, $f^{(2n-1)}(0)=-1$, for all $n\in\N$.
Therefore, for any fixed $m\in\N$, we obtain from (\ref{Muntz-type formula - Re s<-1 - k=1}) and (\ref{Muntz-type formula - Re s<-1 - k=2}) the integral representations
\begin{align}
\zeta(s)\Gamma(s)=\int_{0}^{\infty}
\left(\frac{1}{e^x-1}-\frac{1}{x}+\frac{1}{2}+\sum_{n=1}^{m}
\frac{\zeta(1-2n)}{(2n-1)!}\,x^{2n-1}\right)x^{s-1}dx
\end{align}
and
\begin{align}
\zeta^2(s)\Gamma(s) =\int_{0}^{\infty}\left(\sum_{n=1}^{\infty}d(n)e^{-nx} +\frac{\ln(x)}{x}-\frac{\gamma}{x}-\frac{1}{4}+\sum_{n=1}^{m}
\frac{\zeta^2(1-2n)}{(2n-1)!}\,x^{2n-1}\right)x^{s-1}dx
\end{align}
valid in the strip $-2m-1<\Real(s)<-2m+1$.
\vspace*{0,5 cm}

\subsection*{Acknowledgements}
The author is deeply grateful to Semyon Yakubovich for fruitful discussions of these topics and his useful suggestions, which rather improved the presentation of this paper.

\newpage
\bibliographystyle{plain}
\bibliography{Bibliography}
\end{document}